%% file: 4K1-free-line-graphs-2015-06-24.tex
\documentclass[12pt]{article}
\usepackage{latexsym}
\usepackage{tikz}
\usepackage{tkz-graph}
\newcommand*\circled[1]{\tikz[baseline=(char.base)]{
            \node[shape=circle,draw,inner sep=2pt] (char) {#1};}}
\def\qed{\hfill $\Box$\vspace{2ex}}

\usetikzlibrary{shapes}

\newtheorem{Theorem}{Theorem}[section]
\newtheorem{Corollary}[Theorem]{Corollary}
\newtheorem{Claim}[Theorem]{Claim}

\newtheorem{Lemma}[Theorem]{Lemma}
\newtheorem{Problem}[Theorem]{Problem}

\def\inst#1{$^{#1}$}

\def\FREE{{\it Free}}

\title{A coloring algorithm for $4K_1$-free line graphs}

\author{
    Dallas J. Fraser\inst{1}
    \and Ang\`ele M. Hamel\inst{1}
    \and Ch\'inh T. Ho\`ang\inst{1}\footnote{phone:(519)884-0710, fax:(519)746-0677, email:choang@wlu.ca}
}
\begin{document}
\maketitle

\begin{center}
{\footnotesize

\inst{1} Department of Physics and Computer Science, Wilfrid
Laurier University, \\Waterloo, Ontario, Canada}

\end{center}

\begin{abstract}
Let $L$ be a set of graphs. \FREE($L$) is the set of graphs that
do not contain any graph in $L$ as an induced subgraph. It is
known that if $L$ is a set of four-vertex graphs, then the
complexity of the coloring problem for \FREE($L$) is known with
three exceptions: $L = \{$claw, $4K_1\}$, $L$ = \{claw, $4K_1$,
co-diamond\}, and $L$ = \{$C_4$, $4K_1$\}. In this paper, we study
the coloring problem for \FREE(claw, $4K_1$). We solve the
coloring problem for a subclass of \FREE(claw, $4K_1$) which
contains the class of $4K_1$-free line graphs. Our result implies
the chromatic index of a graph with no matching of size four can
be computed in polynomial time.

\noindent{\em Keywords}: Graph coloring, $claw$, $K_5-e$, $line$
$graph$
\end{abstract}

\section{Introduction}\label{sec:intro}

Graph coloring is one of the most important problems in graph
theory and computer science. Determining the chromatic number of a
graph is a NP-hard problem. However, for some graph families the
problem can be solved in polynomial time. Let $L$ be a set of
graphs. Define \FREE($L$) to be the class of graphs that do not
contain any graph in the list $L$ as an induced subgraph. For
example, \FREE($P_4$) is the class of graphs that do not contain a
$P_4$ as an induced subgraph; and \FREE($P_5$, co-$P_5$) is the
class of graphs that do not contain an induced subgraph isomorphic
to a $P_5$ or the complement of a $P_5$. For a single graph $H$,
in \cite{KraKra2001}, it is proved the coloring problem for
\FREE($H$) is polynomial time solvable if $H$ is the $P_4$, or the
join of a $P_3$ and $P_1$, and NP-complete for any other graph
$H$. This result motivates us to consider the problem of coloring
the class \FREE($L$) when $L$ is a family of four-vertex graphs.
As this paper is being written, we discovered that
\cite{LozMal2015} has considered the same problem. We found some
results already discovered in \cite{LozMal2015}, but we also found
a new result which we will present in this paper. To explain this
result, we will need to discuss the background of the problem. Let
VERTEX COLORING be the problem of determining the chromatic number
of a graph. For graphs $G$ and $H$, $G + H$ denotes the disjoint
union of $G$ and $H$. For a graph $G$ and an integer $k$, $kG$
denotes the disjoint union of $k$ copies of $G$. Let $P_n$
(respectively, $C_n, K_n$) denote the chordless path
(respectively, chordless cycle, clique) on $n$ vertices.

Recall the following result in \cite{KraKra2001}:
\begin{Theorem}\label{thm:KraKra}
For a single graph $H$, VERTEX COLORING is polynomial time
solvable for \FREE($H$) if $H$ is the $P_4$, or $P_3 + P_1$, and
NP-Complete otherwise.
\end{Theorem}

In \cite{Rao2007}, the following result is established (see
\cite{Rao2007} for the definition of clique widths.  For the
purposes of this paper, we need only know the fact that if the
clique width of a graph is bounded then so is that of its
complement.)
\begin{Theorem}\label{thm:Rao2007}
VERTEX COLORING is polynomial time solvable for graphs with
bounded clique width.
\end{Theorem}
\begin{figure}\label{fig:4vertex-graphs}
\input{All4VertexGraphs.tex}
\caption{All 4-vertex graphs}
\end{figure}
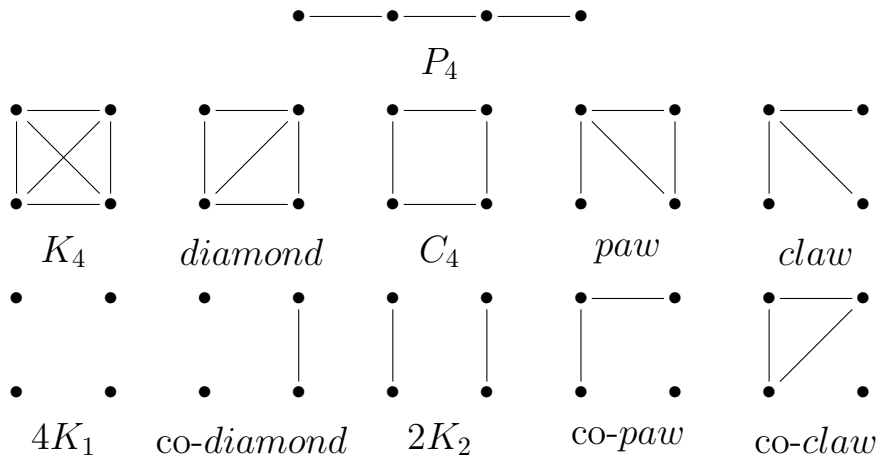
In \cite{BraEng2006}, the authors study the clique widths of
\FREE($F$) where $F$ is a family of four-vertex graphs.
Figure~\ref{fig:4vertex-graphs} shows all 11 graphs on four
vertices with their names. They show there are exactly seven
minimal classes with unbounded clique width. These are:
\newline ${\cal X}_1$ = \FREE(claw, $C_4$, $K_4$, diamond).
\newline ${\cal X}_2$ = \FREE(co-claw, $2K_2$, $4K_1$, co-diamond).
\newline ${\cal X}_3$ = \FREE($C_4$, co-claw, paw, diamond, $K_4$).
\newline ${\cal X}_4$ = \FREE($2K_2$, claw, co-paw, co-diamond, $4K_1$).
\newline ${\cal X}_5$ = \FREE($K_4$, $2K_2$).
\newline ${\cal X}_6$ = \FREE($C_4$, $2K_2$).
\newline ${\cal X}_7$ = \FREE($C_4$, $4K_1$).

Thus, if $F$ is a set of four-vertex graphs and $F \not \subseteq
{\cal X}_i$ ($i = 1, 2, \ldots, 7$), then VERTEX COLORING is
polynomial time solvable for {\it Free}($F$).

VERTEX COLORING is NP-Complete for
\begin{itemize}
 \item ${\cal X}_1$ due to a result in  \cite{KraKra2001} which shows the problem
  is NP-Complete for \FREE(claw, $C_4$) and for \FREE(claw,
$K_4$, diamond);
 \item ${\cal X}_2$ due to Theorem~6 in \cite{Sch2005};
 \item ${\cal X}_3$ due to a remark (Case 1 in Section 4) in \cite{KraKra2001}:
    they show the problem is NP-Complete for \FREE($L$)
    if every graph in $L$ contains a cycle .
\end{itemize}
Thus, VERTEX COLORING is NP-Complete for \FREE($L$) whenever $L
\subseteq {\cal X}_i$, for $i = 1, 2, 3$. Therefore, we only need
to examine the problem for classes ${\cal X}_4$, ${\cal X}_5$,
${\cal X}_6$, ${\cal X}_7$ and their super classes defined by
forbidding induced subgraphs with four vertices. In
\cite{LozMal2015}, a polynomial time algorithm is given for VERTEX
COLORING for the class ${\cal X}_5$. The graphs in ${\cal X}_6$
have a simple structure \cite{BlaHuj1993} that implies an easy
polynomial time algorithm for the coloring problem. Furthermore,
in \cite{HoaLaz2015} a polynomial time algorithm for VERTEX
COLORING for the larger class \FREE($P_5$, co-$P_5$) is given. The
complexity of VERTEX COLORING for the class ${\cal X}_7$ is
unknown. In \cite{LozMal2015}, the authors conjecture that the
coloring problem can be solved in polynomial time for ${\cal
X}_7$.

We are interested in the class ${\cal X}_4$.  Let $H$ be a subset
of \{$2K_2$, claw, co-paw, co-diamond, $4K_1$\}. We examine the
complexity of VERTEX COLORING for \FREE($H$). We may assume $H$
does not contain a co-paw, for otherwise the problem is polynomial
time solvable, by Theorem~\ref{thm:KraKra}. We may assume $H$
contains the claw, for otherwise $H \subset {\cal X}_2$, and so
VERTEX COLORING is NP-Complete.  In \cite{LozMal2015}, a
polynomial time algorithm is given for VERTEX COLORING for class
{\it Free}(claw, $2K_2$).
Thus we have $H \subseteq \{$claw, co-diamond, $4K_1$\}. In
\cite{LozMal2015}, it is proved VERTEX COLORING for \FREE(claw,
co-diamond) is polynomially equivalent to the same problem for the
class \FREE(claw, co-diamond, $4K_1$). Thus, if VERTEX COLORING is
polynomial time solvable for \FREE(claw, $4K_1$), then so is the
same problem for \FREE(claw, co-diamond). These are two
challenging problems. We believe VERTEX COLORING can be solved in
polynomial time for \FREE(claw, $4K_1$). The purpose of this paper
is to solve the problem for a subclass of \FREE(claw, $4K_1$): the
class of $4K_1$-free line graphs.

\begin{figure}
\input{forbidden-linegraphs.tex}

\caption{All minimal forbidden induced subgraphs for line
graphs}\label{fig:linegraph}
\end{figure}
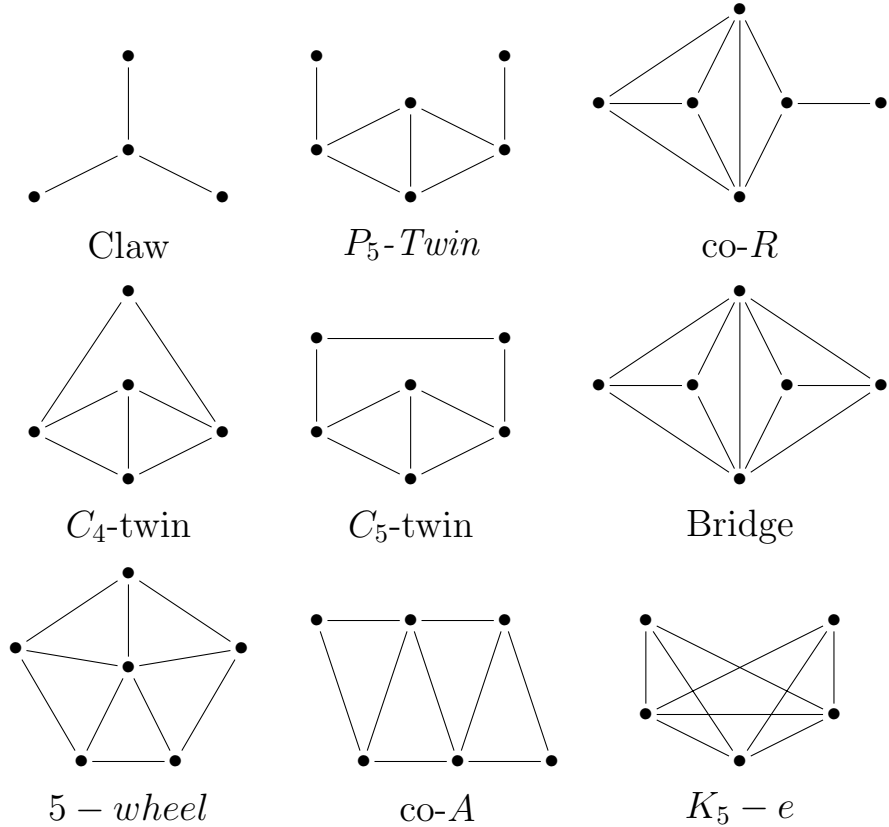

Given a graph $G$, the line graph $L(G)$ of $G$ is defined to be
the graph whose vertices are the edges of $G$, and two vertices of
$L(G)$ are adjacent if their corresponding edges in $G$ are
incident. Line graphs cannot contain a claw. In \cite{Bei1970}, a
characterization of line graphs by forbidden induced subgraphs is
given: A graph is a line graph (of some other graph) if and only
if it does not contain a graph in Figure~\ref{fig:linegraph} as an
induced subgraph.
The purpose of this paper is to prove the following theorem (see
Figure~\ref{fig:linegraph} for the names of the graphs mentioned
in our theorem.)
\begin{Theorem} \label{thm:main}
VERTEX COLORING is polynomial time solvable for \FREE(claw,
$4K_1$, 5-wheel, $C_5$-twin, $P_5$-twin, $K_5 -e$).
\end{Theorem}
An {\it edge coloring} of a graph is an assignment of colors to
its edges such that every edge receives one color and two edges
receive different colors if they are incident.  The {\it chromatic
index} of a graph is the smallest number of colors needed to color
its edges. Let $\Delta(G)$ be the maximum degree of $G$. Then the
chromatic index of a graph $G$ is at least $\Delta(G)$. A classic
theorem of Vizing \cite{Viz1964} states that the chromatic index
of a graph $G$ is at most $\Delta(G) + 1$. However, computing the
chromatic index of a graph is NP-hard \cite{Hol1981}. Let $G$ be a
graph and $L(G)$ be its line graph. Then the chromatic index of
$G$ is equal to the chromatic number of $L(G)$. A matching of $G$
is a set of edges such that no two edges in it are incident; this
matching of $G$ corresponds to a stable set of the line graph
$L(G)$. It follows that our Theorem~\ref{thm:main} implies the
chromatic index can be computed in polynomial time for graphs
without matching of size four.
\begin{Corollary}
There is a polynomial-time algorithm to compute the chromatic
index of  a graph without a matching of size four.
\end{Corollary}
In Section~\ref{sec:background}, we discuss the background results
needed to prove our main theorem. In
Section~\ref{sec:observations}, we establish structural properties
of the graphs in \FREE(claw, $4K_1$, 5-wheel, $C_5$-twin,
$P_5$-twin, $K_5 -e$). In Section~\ref{sec:algorithm}, we prove
Theorem~\ref{thm:main}. And finally, in
Section~\ref{sec:conclusions}, we discuss open problems related to
our work.
\section{Definitions and background}\label{sec:background}
In this section, we discuss the background of our problem.  Let
$G$ be a graph. Then {\it co-$G$} denotes the complement of $G$. A
{\em clique cutset} of $G$ is a set of vertices $S$ where $G[S]$,
the subgraph of $G$ induced by $S$, is a clique whose removal
disconnects $G$. An {\it atom} is a connected graph containing no
clique cutset. Let $\chi(G)$ denote the chromatic number of a
graph $G$. Let $\omega(G)$ denote the number of vertices in a
largest clique of $G$.   Let $\alpha(G)$ denotes the number of
vertices in a largest stable set of $G$.  A {\it hole} is the
graph $C_k$ with $k \geq 5$. For a hole $H$, a vertex $x$ outside
$H$ is a $k$-vertex (for $H$) if $x$ has exactly $k$ neighbours in
$H$. An {\it anti-hole} is the complement of a hole. A hole is
{\it odd} if it has an odd number of vertices. A graph is {\it
Berge} if it contains no odd holes and no odd anti-holes. For sets
$X,Y$ of vertices, we write $X \;\circled{0}\; Y$ to mean there is
no edge between any vertex in $X$ and any vertex in $Y$, and $X
\;\circled{1}\; Y$ to mean there are all edges between $X$ and
$Y$. We let $|G|$ denotes the number of vertices of $G$. By {\it
claw(a,b,c,d)} we denote the claw with vertices $a,b,c,d$ and
edges $ab, ac, ad$. We say that a stable set is {\it good} if it
meets every largest clique of the graph. A set is {\it big} if it
has at least three elements. For a vertex $x$ of $G$, $N(x)$
denote the set of neighbors of $x$, and $d(x)$ is the {\it degree}
of $x$, that is, $d(x) = |N(x)|$.

Consider  the following procedure to decompose a (connected) graph
$G$. If $G$ has a clique cutset $C$, then $G$ can be decomposed
into subgraphs $G_1 = G[V_1]$ and $G_2 = G[V_2]$ where $V = V_1
\cup V_2$ and $C = V_1 \cap V_2$ (recall that $G[X]$ denotes the
subgraph of $G$ induced by $X$ for a subset $X$ of vertices of
$V(G)$). Given optimal colourings of $G_1, G_2$, we can obtain an
optimal colouring of $G$ by identifying the colouring of $C$ in
$G_1$ with that of $C$ in $G_2$. In particular, we have $\chi(G)=
\max(\chi(G_1),\chi(G_2))$. If $G_i$ ($i \in \{1,2\}$) has a
clique cutset, then we can recursively decompose $G_i$ in the same
way.  This decomposition can be represented by a binary tree
$T(G)$ whose root is $G$ and the two children of $G$ are $G_1$ and
$G_2$, which are in turn the roots of subtrees representing the
decompositions of $G_1$ and $G_2$. Each leaf of $T(G)$ corresponds
to an atom of $G$. Algorithmic aspects of the clique cutset
decomposition are studied in \cite{Tar1985} and \cite{Whi1984}. In
particular, the decomposition tree $T(G)$ can be constructed in
$O(nm)$ time such that the total number atoms is at most $n-1$
\cite{Tar1985} (Here, as usual, $n$ and, respectively, $m$, denote
the number of vertices, respectively, edges, of the input graph).
Thus, to color a graph $G$ in polynomial time, we only need to
color its atoms in polynomial time.

Our result relies on known theorems on perfect claw-free graphs,
and we discuss these results now. A graph $G$ is {\it perfect} if
for each induced subgraph $H$ of $G$, we have $\chi(H) =
\omega(H)$. In \cite{ParRav1976}, it is proved claw-free Berge
graphs are perfect. Claw-free perfect graphs can be recognized in
polynomial time \cite{ChvSbi1998}, and they can be optimally
colored in polynomial time \cite{Hsu1981}. We note
\cite{ChuRob2006} proves that a graph is perfect if and only if it
is Berge, solving a long standing conjecture of \cite{Ber1961}.
Perfect graphs can be recognized in polynomial time
\cite{ChuCor2005}, and they can be optimally colored in polynomial
time \cite{GroLov1984}. In \cite{ChvSbi1998}, the following
result, crucial to our algorithm, is established.
\begin{Lemma}\label{lem:BenRebea}
Let $G$ be a connected claw-free graph with $\alpha(G) \geq 3$. If
$G$ contains an odd anti-hole then $G$ contains a $C_5$.
\end{Lemma}
Finally, we note the well known observation that VERTEX COLORING
is polynomial time solvable for graphs $G$ with $\alpha(G) = 2$;
it is sufficient to find a maximum matching in the complement of
$G$.

\section{Structural properties}\label{sec:observations}
In this section, we establish preliminary results needed to prove
Theorem~\ref{thm:main}. For the Claims in this section, we will
assume $G=(V,E)$ is a connected graph in \FREE(claw, $K_5-e$,
5-wheel, $C_5$-twin, $P_5$-twin, $4K_1$). We begin with the
following easy Claim.

\begin{Claim}\label{cla:cl-free}
If $\alpha(G) \leq 4$, then $G$ contains no $C_\ell,\; \ell \geq
8$.
\end{Claim}
{\it Proof.}   If $G$ contains $C_\ell,\; \ell \geq 8$, then $G$
contains a $4K_1$. \qed

Let $C$ be a hole with vertices $1,2, \ldots, k$ (in the cyclic
order) with $k \geq 5$ (with respect to $C$, the vertices $i$ will
be taken modulo $k$). Define $Y_i$ be the set of vertices with
neighbors $i,i+1,i+2,i+3$. Let $Z_i$ be the set of vertices with
neighbors $i,i+1,i+3,i+4$.

\begin{Claim}\label{cla:c7-bounded}
 If $G$ contains a $C_7$ , then $|G|
\leq 21$.
\end{Claim}
\noindent{\it Proof}. Suppose $G$ contains a $C_7$.  Then $G$ has
no $k$-vertex for $k \in \{0, 1, 2\}$ since $G$ is $4K_1$-free,
has no $k$-vertex for $k \in \{5, 6,7\}$ since $G$ is claw-free,
and $G$ has no 3-vertex, or else $G$ contains a claw or a
$P_5$-twin. Thus, only $4$-vertices may exist and they are one of
the two types $Y_i$s, $Z_i$s defined above.
Let $y_1$ and $y_2$ be two vertices in $ Y_i,\; y_1 \neq y_2$. If
$y_1y_2 \not \in E$ then there is a $claw (i+3, i+4, y_1, y_2)$
and if $y_1y_2 \in E$ then there is a $K_5 - e$ induced by
$\{i,i+1,i+2,y_1,y_2\}$; so we have $|Y_i| \leq 1$. Let $z_1$ and
$z_2$ be two vertices in $ Z_i,\; z_1 \neq z_2$. If $z_1z_2 \not
\in E$ then there is a $claw (i+4, i+5, z_1, z_2)$ and if $z_1z_2
\in E$ then there is a $C_5$-twin induced by $\{ i, i+6, i+5, i+4,
z_1, z_2\}$. It follows that $G$ contains at most $21$ vertices.
\qed

\begin{Claim}\label{cla:c5-kvertex}
If $G$ contains a $C_5$ then $G$ has no $k$-vertex for $ k \in
\{1,3,5\}$.
\end{Claim}
\noindent {\it Proof.} Suppose $G$ contains a $C_5$.  Then $G$ has
no $k$-vertex for $k = 1$ for otherwise $G$ contains a claw, or
for $k = 3$ for otherwise $G$ contains a claw or a $C_5$-twin, or
for $k = 5$ for otherwise $G$ contains a 5-wheel. \qed

For the all the Claims below, we will assume $G$ contains a $C_5$.
Let the $0$-vertex set be denoted by $R$, let $X_i$ be the set of
$2$-vertices with neighbours $i, i+1$ and let $Y_i$ be the set of
$4$-vertices with neighbors $i, i+1, i+2, i+3$. Let $X$ denote the
set of all $2$-vertices and $Y$ denote the set of all
$4$-vertices.  Every vertex of $G - C_5$ belongs to $X \cup Y \cup
R$.

\begin{Claim}\label{cla:c5-4vertex-bounded}
We have $|Y_i| \leq 1$ for all $i$.
\end{Claim}
\noindent {\it Proof.} Let $y_1$ and $y_2$ be two vertices in $
Y_i,\; y_1 \neq y_2$. If $y_1y_2 \not \in E$ then there is a $claw
(i,i-1, y_1,y_2)$. If $y_1y_2 \in E$ then there is a $K_5 - e $
induced by  $\{i,i+1,i+2,y_1,y_2\}$.  \qed

From Claim~\ref{cla:c5-4vertex-bounded}, we have
 $| Y | \leq 5$.
\begin{Claim}\label{cla:yi-cojoin-yi1}
We have $Y_i \;\circled{0}\; Y_{i+1}$ for all $i$.
\end{Claim}
\noindent {\it Proof.} Let $y_1$ be the vertex from  $Y_i$ and
$y_2$ be the vertex from $Y_{i+1}$. If $y_1y_2 \in E$ then there
is a $K_5 - e$ induced by $ \{ i+1, i+2, i+3, y_1, y_2 \}$. \qed

For the following two claims, we will let $x_i$ (respectively,
$y_i$) denote an arbitrary vertex in $X_i$ (respectively, $Y_i$)
for all $i$.

\begin{Claim}\label{cla:xi-joins-yi}
$X_i \; \circled{1} \; Y_i \cup Y_{i+3}$.
\end{Claim}
\noindent {\it Proof.} If $x_iy_i \not \in E$,  then there is a
claw($i, i+4, x_i, y_i$). If $x_iy_{i+3} \not \in E$, then there
is a claw($i+1, i+2, x_i, y_{i+3}$). \qed

\begin{Claim}\label{cla:one-y}
If $X_i \not= \emptyset$, and both $y_i$ and $y_{i+3}$ exist, then
$y_i y_{i+3} \in E$.
\end{Claim}
\noindent {\it Proof.} If $y_i y_{i+3} \not\in E$, then by
Claim~\ref{cla:xi-joins-yi}, the set $ \{x_i, i, i+1, y_i, y_{i+3}
\}$ contains a $K_5 - e$. \qed

\begin{Claim}\label{cla:xi-cojoins-yi}
$X_i \; \circled{0} \; Y_{i+1} \cup Y_{i+2} \cup Y_{i+4}$.
\end{Claim}
\noindent {\it Proof.} If $x_iy_{i+1} \in E$, then there is a
claw($y_{i+1}, x_i, i+2, i+4)$. If $x_iy_{i+4} \in E$, then there
is a claw($y_{i+4}, i+4, i+2, x_i$). If $x_iy_{i+2} \in E$, then
there is a claw($y_{i+2}, i+2, i+4, x_i$). \qed

\begin{Claim}\label{cla:R-cojoin-Y}
We have $R \;\circled{0}\; Y$.
\end{Claim}
\noindent {\it Proof.} If some $y \in Y$ is adjacent to some $r
\in R$,  then there is a claw induced by $y,r$ and some two
neighbours $a, b$ of $y$ in the $C_5$ with $ab \not\in E$. \qed

\begin{Claim}\label{cla:x-not-empty}
If $R \not = \emptyset$, then $X \not = \emptyset$. \qed
\end{Claim}
\noindent {\it Proof}. Since $G$ is connected by the assumption of
this section, there is a path connecting some vertex in $R$ to
some vertex in the $C_5$. By Claim~\ref{cla:R-cojoin-Y}, this path
must contain some vertex in $X$. \qed

\begin{Claim}\label{cla:xi-clique}
Each $X_i$ forms a clique for all $i$.
\end{Claim}
\noindent {\it Proof.} Let $v_1$ and $v_2$ be two vertices in $
X_i,\; v_1 \neq v_2$. If $v_1v_2 \not \in E$ then there is a $claw
(i,i-1, v_1, v_2)$. So $v_1v_2 \in E$ and $X_i$ forms a clique.
\qed

\begin{Claim}\label{cla:r-clique}
The set $R$ induces a clique.
\end{Claim}
\noindent {\it Proof.} Let $r_1$ and $r_2$ be vertices in $R$. If
$r_1r_2 \not \in E$ then the set $\{0, 2, r_1, r_2\}$ induces a
$4K_1$. \qed

\begin{Claim}\label{cla:r-join-xi}
$R \; \circled{1} \; X_i$ for all $i$.
\end{Claim}
\noindent {\it Proof.} Let $x$ be a vertex in $ X_i$ and $r$ be a
vertex in $ R$. If $rx \not \in E$ then there is a $4K_1$ induced
by $\{ r, x, i-1, i+2\}$. \qed

\begin{Claim}\label{cla:r-limits-xi}
If $R \neq \emptyset$ then $|X_i| \leq 2$ for all $i$.
\end{Claim}
\noindent {\it Proof.} Suppose $|X_i| \geq 3$ and $R \neq
\emptyset$. Let $x_1,\; x_2,\;$ and $x_3$ be three distinct
vertices from $X_i$.  By Claim \ref{cla:xi-clique}, the vertices
$x_1,\; x_2,\;$ and $x_3$ form a clique. By Claim
\ref{cla:r-join-xi} there is a $K_5 - e$ induced by $\{r, x_1,
x_2, x_3, i\}$ for a vertex $r \in R$. \qed

\begin{Claim}\label{cla:2xi-cojoins-xj}
A vertex in $ X_i$ cannot have two neighbors in $X_j$  for any two
distinct $i$ and $ j$.
\end{Claim}
\noindent {\it Proof.}  Suppose some vertex $x_i \in X_i $ is
adjacent to two vertices $a,b \in X_j$. We may assume $j=i+1$, or
$j=i+2$. Now, there is a $P_5$-twin induced by $\{i, x_i, a, b,
j+1, j+2 \}$ if $j=i+1$,  or $\{j-1, x_i, a, b, j+1, j+2\}$ if $j
= i+2$. \qed

\begin{Claim}\label{cla:xi-limits-r}
If $X \neq \emptyset$ then $|R| \leq 2$ or $X$ is a clique cutset
of $G$.
\end{Claim}
\noindent {\it Proof.}
Suppose $X \neq \emptyset$ and $|R| \geq 3$. By Claim
\ref{cla:r-join-xi}, $R$ $\circled{1}$ $X$. By
Claim~\ref{cla:R-cojoin-Y}, $X$ is a cutset separating $R$ from
the $C_5$. We may assume $X$ contains two non-adjacent vertices
$v_1, v_2$, for otherwise $X$ is a clique cutset and we are done.
But now, by Claim~\ref{cla:r-clique}, any three vertices in $R$
together with $v_1, v_2$ induce a $K_5 - e$ in $G$. \qed

\begin{Claim}\label{cla:r-bounded}
If $R \neq \emptyset$ then $ |G| \leq 22$ or $G$ contains a clique
cutset.
\end{Claim}
\noindent {\it Proof.} If $|R| \geq 3$ then by Claims
\ref{cla:xi-limits-r} and \ref{cla:x-not-empty}, $X$  is a clique
cutset of $G$.  So we have $|R| \leq 2$. By Claim
\ref{cla:r-limits-xi}, we have $|X_i| \leq 2$ for $i \in
{0,1,2,3,4}$. By Claim \ref{cla:c5-4vertex-bounded}, we have $|Y|
\leq 5$. Then $|G| = |R| + |X| + |Y| + |C_5|  \leq 2 + 10 + 5 + 5
= 22$. \qed

\begin{Claim}\label{cla:xi-neighbor-adjacent}
A vertex in $X_i$ cannot have two non-adjacent neighbors in $X -  X_i$.
\end{Claim}
\noindent {\it Proof.} Suppose some vertex $x_i \in X_i$ have
non-adjacent neighbors $v_1$ and $v_2$ in $X - X_i$. By Claim
\ref{cla:xi-clique}, we have  $v_1 \in X_j,\; v_2 \in X_k,\; j
\neq k$. If $j=i-1$ and $k=i+1$ then there is a $P_5$-twin induced
by $\{v_2, x_i, v_1, i, j, j-1\}$. Now, we may assume $k \not \in
\{i-1, i+1\}$,  it follows there is a claw($x_i, v_2, v_1, i+1) $.
\qed

\begin{Claim}\label{cla:three-consecutive-xi}
Suppose $X_i, X_{i+1}, X_{i+2}$ are each non-empty for some $i$.
If $|X_j| \geq 3$ then $|X_k| = |X_\ell| = 1,$ for $ \{j,k,\ell \}
= \{i,i+1,i+2\}$.
\end{Claim}
\noindent {\it Proof.} Suppose $|X_j| \geq 3$. Suppose some vertex
$x_k \in X_k$ is non-adjacent to some vertex $x_\ell \in X_\ell$.
By Claim~\ref{cla:2xi-cojoins-xj}, there is a vertex $x_j \in X_j$
that is non-adjacent to both $x_k$ and $x_\ell$. But now $\{i+4,
x_j, x_k, x_\ell \}$ induces a $4K_1$. Thus, we have $X_k \;
\circled{1} \; X_\ell$. It follows from
Claim~\ref{cla:2xi-cojoins-xj} that $|X_k | = | X_\ell | = 1$.
\qed

\begin{Claim}\label{cla:coloring-xi-xi2-xi3}
Suppose $ R =  Y =  X_{i+2} = X_{i+4} = \emptyset$  for some $i$.
Then $\chi(G) = \omega(G)$, and an optimal colouring of $G$ can be
found in polynomial time.
\end{Claim}
\noindent {\it Proof.} By induction on the number of vertices. Let
$X_k$ be the set with largest cardinality  among the three sets
$X_i, X_{i+1}, X_{i+3}$. If $|X_k| = 1$, then $\omega(G) = 3, |V|
\leq 8$, and a 3-coloring of $G$ can be trivially found.
Similarly, if $|X_k| = 2$, then $\omega(G) = 4$, and  a 4-coloring
of $G$ can be found. Now, we may assume $|X_k| \geq 3$. From
Claims~\ref{cla:2xi-cojoins-xj}
and~\ref{cla:xi-neighbor-adjacent}, there is a stable set $S$
containing a vertex in $X_k$ and some vertices in $X_j \cup
X_{\ell}$ with $k \not\in \{j, \ell\}$ that is good, ie.,  meets
all maximum cliques of $G$. Now, we can recursively color $G-S$
with $\omega(G) -1$ colors and then give $S$ a new color. \qed

In proving the main theorem of this paper, we will reduce the
problem to list coloring a restricted class of graphs. We will now
define some new notions. Given a graph $G$ and a list of colors
$L(v)$ for each vertex of $v$, an $L$-coloring of $G$ is a
(proper) coloring such that each vertex is assigned a color from
its list.

\begin{Lemma}\label{lem:special-graph}
Let $G=(V,E)$ be a graph whose vertices can be partitioned into
three cliques $Q_1, Q_2, Q_3$ such that
\begin{description}
 \item[(a)] each vertex in $Q_i$ is adjacent to at most one vertex in
$Q_j$ for  all $i,j$ with $i \not = j$.
 \item[(b)] if a vertex in $Q_i$ is adjacent to vertices $b \in Q_j, c
 \in Q_k$, then $bc \in E$ for distinct $i,j,k$.
\item[(c)] for each $Q_i$, there is a list $L_i$ of colors such
that
\begin{description}
 \item[(i)] all vertices $v \in Q_i$ have the same list $L(v) =
L_i$,
 \item[(ii)] $|L_i| \geq |Q_i|$,
 \item[(iii)] each $L_i$ contains a color $d_i$ such that the three colors
 $d_1, d_2, d_3$ are all distinct.
 \item[(iv)] no color appears in all three $L_i$.
\end{description}
\end{description}
Then $G$ admits an $L$-coloring where $L = L_1 \cup L_2 \cup L_3$.
\end{Lemma}
\noindent {\it Proof.}  Let $|L|$ denotes the number of colors in
$L$. We prove the Lemma by induction on $|L|$. Below, $i,j,k$ are
distinct indices taken from $\{1,2,3\}$.

If some clique $Q_i$ is empty, then it is easy to see the Lemma
holds. If $|Q_1| = |Q_2| =|Q_3| = 1$, then color the only vertex
in $Q_i$ with color $d_i$, and we are done. Thus, some $Q_i$ has
at least two vertices.

Suppose some $Q_i$ has exactly one vertex $v$, and  $d_i \not\in
L_j \cup L_k$. In this case,   color $v$ with color $d_i$, remove
$v$ from $G$. By induction, $G-v$ admits an $L$-coloring, and we
are done.

Suppose $Q_i$ has exactly one vertex, $d_i \in L_j$, and $Q_j$ has
at least two vertices. Let $v$ be the only vertex in $Q_i$. Take a
vertex $b \in Q_j$ which is not adjacent to $v$; color $v$ and $b$
with color $d_i$, remove $d_i$ from $L_j$ (note that $d_j$ remains
in $ L_j$). By induction $G-\{v,b\}$ admits an $L$-coloring, and
we are done.

Now, suppose  $Q_i$ has exactly one vertex $v$. From the above, we
may assume $d_i \in Q_j$, and $|Q_j| = 1$. It follows that $|Q_k|
\geq 2$ and thus, $d_i \not\in L_k$. Note that we have $\{d_i, d_j
\} \subseteq L_j$. Color $v$ with color $d_i$, remove color $d_i$
from $L_j$. By induction $G-v$ admits an $L$-coloring, and we are
done.

We may now assume each $Q_i$ has at least two vertices. Suppose
there is a color $c \not\in \{d_1, d_2, d_3\}$. If $c$ appears in
only one $Q_i$, then assign to a vertex $v \in Q_i$ the color $c$,
remove $c$ from $L_i$, and by induction $G-v$ admits an
$L$-coloring and we are done. Now, assume $c \in Q_i$ and $c \in
Q_j$. Since each of $Q_i$ and $Q_j$ has at least two vertices,
there are non adjacent vertices (by condition (a)) $a \in Q_i$, $b
\in Q_j$. Assign to $a$ and $b$ color $c$, remove color $c$ from
$L_i$ and $L_j$. By induction, $G-\{a,b\}$ admits an $L$-coloring,
and we are done.

So, we may assume every color belongs to $\{d_1, d_2, d_3\}$. It
follows that  $2 \leq |Q_i| \leq 3$ for any $Q_i$. By condition
(iv), we have $|Q_i| < 3$ for any $i$. It follows that $|Q_1| =
|Q_2| =|Q_3| = 2$. It is now straightforward to verify that $G$
admits an $L$-coloring.
\qed

\begin{Lemma}\label{lem:k-colorable}
Let $G$ be a graph in \FREE(claw, $4K_1$, 5-wheel, $K_5-e$,
$P_5$-twin, $C_5$-twin) with a $C_5$. Suppose $ R =  X_{i+2} =
X_{i+4} = \emptyset$ for some $i$. If $\omega(G) \geq 5$ then
$\chi(G) = \omega(G)$ and an optimal coloring of $G$ can be found
in polynomial time
\end{Lemma}
\noindent {\it Proof.} Suppose $\omega(G) = 5$. For simplicity, we
may assume $i=1$, i.e., the sets $X_3 $ and $X_5$ are empty. Each
of the sets $X_i$ has at most three vertices for $i=1,2,4$ (since
$X_i \cup \{i, i+1\}$ induces a clique by
Claim~\ref{cla:xi-clique}  ). For each vertex $i$ in the $C_5$,
color vertex $i$ with color $c(i) = i$. Color the vertex $y_i$ (if
it exists) with color $c(y_i) = i -1 $ with the subscript $i$
taken modulo 5. Now, we need to extend this coloring to the set
$X$. For $X_i$, we will construct a list $L_i$ of feasible colors
which are compatible with the already colored vertices in $C_5
\cup Y$. Each vertex in $X_i$ will have the list $L_i$. Consider a
non-empty set $X_i$ ($i \in \{1, 2, 4\}$). If $| X_i | = 1$, then
set $L_i = \{c(i+3)\}$ (which is the color of vertex $i+3$). If
$|X_i | = 2$, then the vertex $y_{i}$ or the vertex $y_{i+3}$ (or
both) does not exist; for otherwise, by Claim~\ref{cla:one-y}, the
induced graph $G[X_i \cup \{i, i+1, y_i, y_{i+3} \}]$ contains a
$K_6$, contradicting the assumption that $\omega(G) = 5$. Let $j$
be the subscript such that $y_j$ ($j \in \{i, i+3\}$) does not
exist. Let $L_i = \{c(i +3), c(y_j)\}$. If $|X_i| = 3$, then both
vertices $y_{i}, y_{i+3}$ do not exist for the same reason above.
Let $L_i = \{c(i +2), c(i+3), c(i+4)\}$. A coloring of the
vertices of $X$ using colors from  the lists $L_i$ does not
conflict with the already assigned coloring of $C_5 \cup Y$. By
Lemma~\ref{lem:special-graph}, $X$ admits an $L$-coloring. So we
have $\chi(G) = 5$.

Now, we may assume $\omega(G) \geq 6$. We will find a good stable
set (recall the definition of good stable sets in
Section~\ref{sec:background}). Note that if a maximum clique $C$
of $G$ contains some vertex of $X_i$, then $C$ contains all
vertices of $X_i$. Also, since $\omega(G) \geq 6$,  a maximum
clique of $G$ does not contain a set $X_i$ with $|X_i| = 1$.  Now,
since $\omega(G) \geq 6$, some set $X_i$ must have at least three
vertices. By Claims~\ref{cla:2xi-cojoins-xj}
and~\ref{cla:xi-neighbor-adjacent}, there is a stable set
containing a vertex in every $X_j$ ($j \in \{1,2,4\}$) with $|X_j|
\geq 2$. Such a set $S$ is the desired stable set. Now, we have
$\omega(G - S) = \omega(G)  -1 $. We recursively color $G-S$ with
$\omega(G) -1 $ colors, and then give vertices in $S$ a new color,
and we are done.
\qed

\begin{figure}\label{fig:petersen}
\input{Petersen.tex}
\caption{The Petersen graph (b) and its complement (a).}
\end{figure}
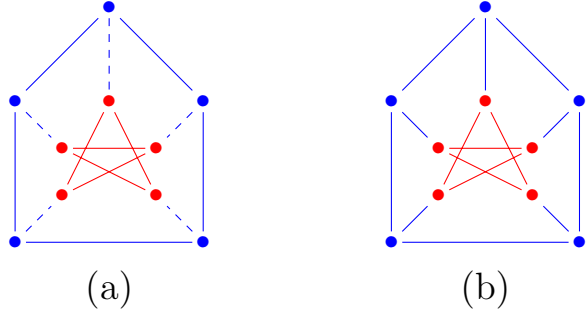

Note that the statement of the Lemma~\ref{lem:k-colorable} is
false for $\omega(G) = 4$.  There are graphs $G$ in \FREE(claw,
$4K_1$, 5-wheel, $K_5-e$, $P_5$-twin, $C_5$-twin) with $\omega(G)
=4$ and $\chi(G) = 5$. The graph in Figure~\ref{fig:petersen} (a)
is such a graph. This is the graph with a $C_5$ (indicated by the
outer $C_5$) and all five $y_i$ vertices (indicated by the inner
$C_5$), there are all edges between the outer and inner $C_5$s
except for the non-edges denoted by the dotted lines. It is very
interesting to note that this graph is the complement of the
Petersen graph.

\section{Coloring algorithm}\label{sec:algorithm}

In the section, we prove Theorem~\ref{thm:main}. Let $G$ be a
graph satisfying the hypothesis of the theorem.  From the
discussion in Section~\ref{sec:background}, we may assume $G$ is
connected, is an atom, ie., $G$ contains no clique cutset, and has
$\alpha(G) \geq 3$. Furthermore, if $G$ is perfect, then we are
done by \cite{Hsu1981} or \cite{GroLov1984}. If $G$ contains an
odd anti-hole, then by Lemma~\ref{lem:BenRebea}, $G$ contains a
$C_5$. So, we may assume $G$ contains an odd hole $H$. Since $G$
is $4K_1$-free, $H$ has 5 or 7 vertices. If $H$ has 7 vertices,
then by Claim~\ref{cla:c7-bounded}, $|G|$ is a constant and we are
done. So, $H$ is a $C_5$. From this $C_5$, define the sets $X,Y,R$
as before, and we may rely on the Claims in
Section~\ref{sec:observations}. We may assume $|G|$ is not a
constant. By Claim~\ref{cla:r-bounded}, we have $R = \emptyset$.

We may assume $\omega(G) \geq 14$; for otherwise,  Ramsey's
theorem \cite{Ram1929} shows that $|G|$ is a constant (both
$\omega(G)$ and $\alpha(G)$ are constants, so $|G|$ is a
constant). If there is a vertex $v$ with degree $d(v) \leq 13$,
then we recursively color $G-v$ optimally, and then give $v$ a
color not appearing in $N(v)$; such a color exists because
$\chi(G-v) \geq \omega(G-v) = \omega(G) \geq 14$. So, we may
assume every vertex of $G$ has degree at least 14.

Suppose some non-empty $X_i$ has $|X_i| \leq 2$. Then for any
vertex $x_i \in X_i$, we have $d(x_i) = |C_5 \cap N(x_i) | + |Y
\cap N(x_i) | + |N(x_i) \cap X| \leq 2 + 5 + 5=  12$ by
Claims~\ref{cla:2xi-cojoins-xj} and~\ref{cla:c5-4vertex-bounded}.
Thus, if $X_i$ is non-empty then it is big (has at least three
vertices.)

For a vertex $i$ in the $C_5$, we may assume $X_i$ or $X_{i-1}$ is
big;  for otherwise, $d(i) = |C_5 \cap N(i)| + |X_i| + |X_{i+1}| +
|Y \cap N(i)| \leq 2 + 2 + 2 + 5 =  11$ by
Claim~\ref{cla:c5-4vertex-bounded}. Thus, at least three sets
$X_i$'s must be big. It follows from
Claim~\ref{cla:three-consecutive-xi} that precisely three sets
$X_i$ are big, and they are not consecutive. By
Lemma~\ref{lem:k-colorable}, we can color $G$ with $\omega(G)$
colors in polynomial time. \qed

\section{Conclusions and open problems}\label{sec:conclusions}
Let $L$ be a family of four-vertex graphs. As mentioned in
Section~\ref{sec:intro}, the complexity of VERTEX COLORING is
known for the class \FREE($L$) with three exceptions: $L$ =
\{claw, $4K_1\}$, $L$ = \{claw, $4K_1$, co-diamond \}, and $L$ =
\{$C_4$, $4K_1$\}. We believe each of the three problems is
polynomial time solvable. In this paper, we studied the problem
for \FREE(claw, $4K_1$). We solved the coloring problem for a
subclass, the class of $4K_1$-free line graphs. Our result implies
the chromatic index of a graph with no matching of size 4 can be
computed in polynomial time. This is an interesting result in its
own right. We conclude our paper with the two open problems below.

\begin{Problem}
For each fixed $k$, is there a polynomial time algorithm to
compute the chromatic index of a graph without a matching of size
$k$?
\end{Problem}

\begin{Problem}
For each fixed $k$, is there a polynomial time algorithm to solve
VERTEX COLORING for the class \FREE(claw, $kK_1$)?
\end{Problem}

\begin{center}
{\bf Acknowledgement}
\end{center}
This work was supported by the Canadian Tri-Council Research
Support Fund. The authors A.M.H. and C.T.H. were each supported by
individual NSERC Discovery Grants. Author D.J.F was supported by
an NSERC Undergraduate Student Research Award (USRA).

\end{document}

%% file: All4VertexGraphs.tex
\begin{center}
\begin{tikzpicture} [scale = 1.25]
\tikzstyle{every node}=[font=\small]

\newcommand{\size}{1}

\newcommand{\p4}{1}{
    \path (\size * 4, 0) coordinate (g1);
    \path (g1) +(-\size, 0) node (g1_1){};
    \path (g1) +(0, 0) node (g1_2){};
    \path (g1) +(\size, 0) node (g1_3){};
    \path (g1) +(2 * \size, 0) node (g1_4){};
    \foreach \Point in {(g1_1), (g1_2), (g1_3), (g1_4)}{
        \node at \Point{\textbullet};
    }
    \draw   (g1_1) -- (g1_2)
            (g1_2) -- (g1_3)
            (g1_3) -- (g1_4);
    \path (g1) ++(\size  / 2,-\size / 2) node[draw=none,fill=none] { {\large $P_4$}};
}

\newcommand{\kfour}{2}{
    \path (0, - \size * 2) coordinate (g2);
    \path (g2) +(0, 0) node (g2_1){};
    \path (g2) +(0, \size) node (g2_2){};
    \path (g2) +(\size, \size) node (g2_3){};
    \path (g2) +(\size, 0) node (g2_4){};
    \foreach \Point in {(g2_1), (g2_2), (g2_3), (g2_4)}{
        \node at \Point{\textbullet};
    }
    \draw   (g2_1) -- (g2_2)
            (g2_1) -- (g2_3)
            (g2_1) -- (g2_4)
            (g2_2) -- (g2_3)
            (g2_2) -- (g2_4)
            (g2_3) -- (g2_4);
    \path (g2) ++(\size  / 2,-\size / 2) node[draw=none,fill=none] { {\large $K_4$}};
}

\newcommand{\diam}{3}{
    \path (\size * 2, - \size * 2) coordinate (g3);
    \path (g3) +(0, 0) node (g3_1){};
    \path (g3) +(0, \size) node (g3_2){};
    \path (g3) +(\size, \size) node (g3_3){};
    \path (g3) +(\size, 0) node (g3_4){};
    \foreach \Point in {(g3_1), (g3_2), (g3_3), (g3_4)}{
        \node at \Point{\textbullet};
    }
    \draw   (g3_1) -- (g3_2)
            (g3_1) -- (g3_3)
            (g3_1) -- (g3_4)
            (g3_2) -- (g3_3)
            (g3_3) -- (g3_4);
    \path (g3) ++(\size  / 2,-\size / 2) node[draw=none,fill=none] { {\large $diamond$}};
}

\newcommand{\cfour}{4}{
    \path (\size * 4, - \size * 2) coordinate (g4);
    \path (g4) +(0, 0) node (g4_1){};
    \path (g4) +(0, \size) node (g4_2){};
    \path (g4) +(\size, \size) node (g4_3){};
    \path (g4) +(\size, 0) node (g4_4){};
    \foreach \Point in {(g4_1), (g4_2), (g4_3), (g4_4)}{
        \node at \Point{\textbullet};
    }
    \draw   (g4_1) -- (g4_2)
            (g4_1) -- (g4_4)
            (g4_2) -- (g4_3)
            (g4_3) -- (g4_4);
    \path (g4) ++(\size  / 2,-\size / 2) node[draw=none,fill=none] { {\large $C_4$}};
}

\newcommand{\paw}{5}{
    \path (\size * 6, - \size * 2) coordinate (g5);
    \path (g5) +(0, 0) node (g5_1){};
    \path (g5) +(0, \size) node (g5_2){};
    \path (g5) +(\size, \size) node (g5_3){};
    \path (g5) +(\size, 0) node (g5_4){};
    \foreach \Point in {(g5_1), (g5_2), (g5_3), (g5_4)}{
        \node at \Point{\textbullet};
    }
    \draw   (g5_1) -- (g5_2)
            (g5_2) -- (g5_3)
            (g5_2) -- (g5_4)
            (g5_3) -- (g5_4);
    \path (g5) ++(\size  / 2,-\size / 2) node[draw=none,fill=none] { {\large $paw$}};
}

\newcommand{\claw}{6}{
    \path (\size * 8, - \size * 2) coordinate (g6);
    \path (g6) +(0, 0) node (g6_1){};
    \path (g6) +(0, \size) node (g6_2){};
    \path (g6) +(\size, \size) node (g6_3){};
    \path (g6) +(\size, 0) node (g6_4){};
    \foreach \Point in {(g6_1), (g6_2), (g6_3), (g6_4)}{
        \node at \Point{\textbullet};
    }
    \draw   (g6_1) -- (g6_2)
            (g6_2) -- (g6_4)
            (g6_2) -- (g6_3);
    \path (g6) ++(\size  / 2,-\size / 2) node[draw=none,fill=none] { {\large $claw$}};
}

\newcommand{\cokfour}{7}{
    \path (0, - \size * 4) coordinate (g7);
    \path (g7) +(0, 0) node (g7_1){};
    \path (g7) +(0, \size) node (g7_2){};
    \path (g7) +(\size, \size) node (g7_3){};
    \path (g7) +(\size, 0) node (g7_4){};
    \foreach \Point in {(g7_1), (g7_2), (g7_3), (g7_4)}{
        \node at \Point{\textbullet};
    }

    \path (g7) ++(\size  / 2,-\size / 2) node[draw=none,fill=none] { {\large $4K_1$}};
}

\newcommand{\codiamond}{8}{
    \path (\size * 2, - \size * 4) coordinate (g8);
    \path (g8) +(0, 0) node (g8_1){};
    \path (g8) +(0, \size) node (g8_2){};
    \path (g8) +(\size, \size) node (g8_3){};
    \path (g8) +(\size, 0) node (g8_4){};
    \foreach \Point in {(g8_1), (g8_2), (g8_3), (g8_4)}{
        \node at \Point{\textbullet};
    }
    \draw   (g8_3) -- (g8_4);
    \path (g8) ++(\size  / 2,-\size / 2) node[draw=none,fill=none] { {\large co-$diamond$}};
}

\newcommand{\cocfour}{9}{
    \path (\size * 4, - \size * 4) coordinate (g8);
    \path (g8) +(0, 0) node (g8_1){};
    \path (g8) +(0, \size) node (g8_2){};
    \path (g8) +(\size, \size) node (g8_3){};
    \path (g8) +(\size, 0) node (g8_4){};
    \foreach \Point in {(g8_1), (g8_2), (g8_3), (g8_4)}{
        \node at \Point{\textbullet};
    }
    \draw (g8_1) -- (g8_2)
          (g8_3) -- (g8_4);
    \path (g8) ++(\size  / 2,-\size / 2) node[draw=none,fill=none] { {\large $2K_2$}};
}

\newcommand{\copaw}{10}{
    \path (\size * 6, - \size * 4) coordinate (g9);
    \path (g9) +(0, 0) node (g9_1){};
    \path (g9) +(0, \size) node (g9_2){};
    \path (g9) +(\size, \size) node (g9_3){};
    \path (g9) +(\size, 0) node (g9_4){};
    \foreach \Point in {(g9_1), (g9_2), (g9_3), (g9_4)}{
        \node at \Point{\textbullet};
    }
    \draw   (g9_1) -- (g9_2)
            (g9_2) -- (g9_3);
    \path (g9) ++(\size  / 2,-\size / 2) node[draw=none,fill=none] { {\large co-$paw$}};
}

\newcommand{\coclaw}{11}{
    \path (\size * 8, - \size * 4) coordinate (g10);
    \path (g10) +(0, 0) node (g10_1){};
    \path (g10) +(0, \size) node (g10_2){};
    \path (g10) +(\size, \size) node (g10_3){};
    \path (g10) +(\size, 0) node (g10_4){};
    \foreach \Point in {(g10_1), (g10_2), (g10_3), (g10_4)}{
        \node at \Point{\textbullet};
    }
    \draw   (g10_1) -- (g10_2)
            (g10_1) -- (g10_3)
            (g10_2) -- (g10_3);
    \path (g10) ++(\size  / 2,-\size / 2) node[draw=none,fill=none] { {\large co-$claw$}};
}

\end{tikzpicture}
\end{center}

%% file: forbidden-linegraphs.tex
\begin{center}
\begin{tikzpicture} [scale = 1.25]
\tikzstyle{every node}=[font=\small]

\newcommand{\size}{1}

\newcommand{\claw}{1}{
  \path (0,0) coordinate (g1);
  \path (g1) +(\size / 2, \size / 2) node (g1_1){}; \path (g1) +(\size / 2, 1.5 * \size) node
  (g1_2){}; \path (g1) +(-\size / 2, 0) node (g1_3){}; \path (g1)
  +(1.5 * \size, 0) node (g1_4){};
  \foreach \Point in {(g1_1),(g1_2),(g1_3),(g1_4)}{
      \node at \Point {\textbullet};
  }
  \draw   (g1_1) -- (g1_2)
          (g1_1) -- (g1_3)
          (g1_1) -- (g1_4);
  \path (g1) ++(\size  / 2,-\size / 2) node[draw=none,fill=none] { {\large Claw}};
}

\newcommand{\g2}{2}{
    \path(3 * \size, 0) coordinate(g2);
    \path (g2) +(\size / 2, 0) node(g2_1){};
    \path (g2) +(- \size / 2, \size / 2) node(g2_2){};
    \path (g2) +(1.5 * \size, \size / 2) node(g2_3){};
    \path (g2) +(\size / 2, \size) node(g2_4){};
    \path (g2) +(- \size / 2, 1.5 * \size ) node(g2_5){};
    \path (g2) +(1.5 * \size, 1.5 * \size) node(g2_6){};

    \foreach \Point in {(g2_1),(g2_2),(g2_3),(g2_4), (g2_5), (g2_6)}{
      \node at \Point {\textbullet};
    }
    \draw   (g2_1) -- (g2_2)
            (g2_1) -- (g2_3)
            (g2_1) -- (g2_4)
            (g2_4) -- (g2_2)
            (g2_4) -- (g2_3)
            (g2_2) -- (g2_5)
            (g2_3) -- (g2_6);
    \path (g2) ++(\size  / 2,-\size / 2) node[draw=none,fill=none] { {\large {\it $P_5$-Twin}}};
}

\newcommand{\g3}{3}{
    \path( 6.5 * \size, 0) coordinate(g3);
    \path(g3) +(\size / 2, 0) node(g3_1){};
    \path(g3) +(0, \size) node(g3_2){};
    \path(g3) +(\size, \size) node(g3_3){};
    \path(g3) +(\size / 2, 2 * \size) node(g3_4){};
    \path(g3) +(-1 * \size, \size) node(g3_5){};
    \path(g3) +(2 * \size, \size) node(g3_6){};

    \foreach \Point in {(g3_1),(g3_2),(g3_3),(g3_4), (g3_5), (g3_6)}{
      \node at \Point {\textbullet};
    }
    \draw   (g3_1) -- (g3_2)
            (g3_1) -- (g3_3)
            (g3_1) -- (g3_4)
            (g3_4) -- (g3_2)
            (g3_4) -- (g3_3)
            (g3_1) -- (g3_5)
            (g3_2) -- (g3_5)
            (g3_4) -- (g3_5)
            (g3_3) -- (g3_6);
    \path (g3) ++(\size  / 2,-\size / 2) node[draw=none,fill=none] { {\large co-$R $}};

}

\newcommand{\g4}{4}{
    \path(0, -3 * \size) coordinate(g4);
    \path (g4) +(\size / 2, 0) node(g4_1){};
    \path (g4) +(- \size / 2, \size / 2) node(g4_2){};
    \path (g4) +(1.5 * \size, \size / 2) node(g4_3){};
    \path (g4) +(\size / 2, \size) node(g4_4){};
    \path (g4) +(\size / 2, 2 * \size ) node(g4_5){};

    \foreach \Point in {(g4_1),(g4_2),(g4_3),(g4_4), (g4_5)}{
      \node at \Point {\textbullet};
    }
    \draw   (g4_1) -- (g4_2)
            (g4_1) -- (g4_3)
            (g4_4) -- (g4_2)
            (g4_4) -- (g4_3)
            (g4_2) -- (g4_5)
            (g4_3) -- (g4_5)
            (g4_1) -- (g4_4);
    \path (g4) ++(\size  / 2,-\size / 2) node[draw=none,fill=none] { {\large $C_4$-twin}};
}

\newcommand{\g5}{5}{
    \path(3 * \size, -3 * \size) coordinate(g5);
    \path (g5) +(\size / 2, 0) node(g5_1){};
    \path (g5) +(- \size / 2, \size / 2) node(g5_2){};
    \path (g5) +(1.5 * \size, \size / 2) node(g5_3){};
    \path (g5) +(\size / 2, \size) node(g5_4){};
    \path (g5) +(- \size / 2, 1.5 * \size ) node(g5_5){};
    \path (g5) +(1.5 * \size, 1.5 * \size) node(g5_6){};

    \foreach \Point in {(g5_1),(g5_2),(g5_3),(g5_4), (g5_5), (g5_6)}{
      \node at \Point {\textbullet};
    }
    \draw   (g5_1) -- (g5_2)
            (g5_1) -- (g5_3)
            (g5_4) -- (g5_2)
            (g5_4) -- (g5_3)
            (g5_2) -- (g5_5)
            (g5_3) -- (g5_6)
            (g5_5) -- (g5_6)
            (g5_1) -- (g5_4);
    \path (g5) ++(\size  / 2,-\size / 2) node[draw=none,fill=none] { {\large $C_5$-twin}};
}

\newcommand{\g6}{6}{
    \path( 6.5 * \size, -3 * \size) coordinate(g6);
    \path(g6) +(\size / 2, 0) node(g6_1){};
    \path(g6) +(0, \size) node(g6_2){};
    \path(g6) +(\size, \size) node(g6_3){};
    \path(g6) +(\size / 2, 2 * \size) node(g6_4){};
    \path(g6) +(-1 * \size, \size) node(g6_5){};
    \path(g6) +(2 * \size, \size) node(g6_6){};

    \foreach \Point in {(g6_1),(g6_2),(g6_3),(g6_4), (g6_5), (g6_6)}{
      \node at \Point {\textbullet};
    }
    \draw   (g6_1) -- (g6_2)
            (g6_1) -- (g6_3)
            (g6_1) -- (g6_4)
            (g6_4) -- (g6_2)
            (g6_4) -- (g6_3)
            (g6_1) -- (g6_5)
            (g6_2) -- (g6_5)
            (g6_4) -- (g6_5)
            (g6_1) -- (g6_6)
            (g6_3) -- (g6_6)
            (g6_4) -- (g6_6);
    \path (g6) ++(\size  / 2,-\size / 2) node[draw=none,fill=none] { {\large Bridge}};

}

\newcommand{\g7}{7}{
    \path(\size / 2, -6 * \size) coordinate(g7);
    \path(g7) +(0, \size) node(g7_1){};
    \path(g7) +(\size / 2, 0) node(g7_2){};
    \path(g7) +(-\size / 2 , 0) node(g7_3){};
    \path(g7) +(-6 / 5 * \size, 6/5 * \size) node(g7_4){};
    \path(g7) +(0, 2 * \size) node(g7_5){};
    \path(g7) +(6 / 5 * \size, 6/5 * \size) node(g7_6){};
    \foreach \Point in {(g7_1),(g7_2),(g7_3),(g7_4), (g7_5), (g7_6)}{
      \node at \Point {\textbullet};
    }
    \draw   (g7_1) -- (g7_2)
            (g7_1) -- (g7_3)
            (g7_1) -- (g7_4)
            (g7_1) -- (g7_5)
            (g7_1) -- (g7_6)
            (g7_2) -- (g7_3)
            (g7_2) -- (g7_6)
            (g7_3) -- (g7_4)
            (g7_4) -- (g7_5)
            (g7_5) -- (g7_6);
    \path (g7) ++(0,-\size / 2) node[draw=none,fill=none] { {\large $5-wheel$}};

}

\newcommand{\g8}{8}{
    \path(\size * 3.5, -6 * \size) coordinate(g8);
    \path(g8) +(-\size, 1.5 * \size) node(g8_1){};
    \path(g8) +(0, 1.5 *  \size) node(g8_2){};
    \path(g8) +(\size , 1.5 * \size) node(g8_3){};
    \path(g8) +(-\size / 2, 0) node(g8_4){};
    \path(g8) +(\size / 2, 0) node(g8_5){};
    \path(g8) +(1.5 * \size, 0) node(g8_6){};
    \foreach \Point in {(g8_1),(g8_2),(g8_3),(g8_4), (g8_5), (g8_6)}{
      \node at \Point {\textbullet};
    }
    \draw   (g8_1) -- (g8_2)
            (g8_1) -- (g8_4)
            (g8_2) -- (g8_3)
            (g8_2) -- (g8_4)
            (g8_2) -- (g8_5)
            (g8_3) -- (g8_5)
            (g8_3) -- (g8_6)
            (g8_4) -- (g8_5)
            (g8_5) -- (g8_6);
    \path (g8) ++(\size  / 3.5,-\size / 2) node[draw=none,fill=none] { {\large co-$A$}};
}

\newcommand{\g9}{9}{
    \path(\size * 7.0, -6 * \size) coordinate(g9);
    \path(g9) +(0, 0) node(g9_1){};
    \path(g9) +(-\size, \size/2) node(g9_2){};
    \path(g9) +(\size, \size/2) node(g9_3){};
    \path(g9) +(-\size, 1.5 * \size) node(g9_4){};
    \path(g9) +(\size, 1.5 * \size) node(g9_5){};
    \foreach \Point in {(g9_1),(g9_2),(g9_3),(g9_4), (g9_5)}{
      \node at \Point {\textbullet};
    }
    \draw   (g9_1) -- (g9_2)
            (g9_1) -- (g9_3)
            (g9_1) -- (g9_4)
            (g9_1) -- (g9_5)
            (g9_2) -- (g9_3)
            (g9_2) -- (g9_4)
            (g9_2) -- (g9_5)
            (g9_3) -- (g9_4)
            (g9_3) -- (g9_5);
    \path (g9) ++(0,-\size / 2) node[draw=none,fill=none] { {\large $K_5 - e$}};
}

\end{tikzpicture}
\end{center}

%% file: Petersen.tex
\begin{center}
\begin{tikzpicture} [scale = 1.25]
\tikzstyle{every node}=[font=\small]

\newcommand{\size}{1}

\newcommand{\copeterson}{1}{
    \path (0, 0) coordinate (g1);
    \path (g1) +(0, 0) node (g1_1){};
    \path (g1) +(0, 3/2 * \size) node (g1_2){};
    \path (g1) +(\size, 5/2 * \size) node (g1_3){};
    \path (g1) +(2 * \size, 3/2 * \size) node (g1_4){};
    \path (g1) +(2 * \size, 0) node (g1_5){};
    \path (g1) +(1/2 * \size, 1/2 * \size) node (g1_6){};
    \path (g1) +(1/2 * \size, \size) node (g1_7){};
    \path (g1) +(\size, 3/2 * \size) node (g1_8){};
    \path (g1) +(3/2 * \size, \size) node (g1_9){};
    \path (g1) +(3/2 * \size, 1/2 * \size) node (g1_10){};

    \foreach \Point in {(g1_1), (g1_2), (g1_3), (g1_4), (g1_5)}{
        \node at \Point{\textcolor{blue}{\textbullet}};
    }
    \foreach \Point in {(g1_6), (g1_7), (g1_8), (g1_9), (g1_10)}{
        \node at \Point{\textcolor{red}{\textbullet}};
    }
    \draw [blue]
            (g1_1) -- (g1_2)
            (g1_1) -- (g1_5)
            (g1_2) -- (g1_3)
            (g1_3) -- (g1_4)
            (g1_4) -- (g1_5);
    \draw [red]
            (g1_6) -- (g1_8)
            (g1_6) -- (g1_9)
            (g1_7) -- (g1_9)
            (g1_7) -- (g1_10)
            (g1_8) -- (g1_10);
    \draw [dashed, blue]
            (g1_1) -- (g1_6)
            (g1_2) -- (g1_7)
            (g1_3) -- (g1_8)
            (g1_4) -- (g1_9)
            (g1_5) -- (g1_10);
    \path (g1) ++(\size,-\size / 2) node[draw=none,fill=none] { {\large (a)}};

}

\newcommand{\peterson}{2}{
    \path (\size * 4, 0) coordinate (g2);
    \path (g2) +(0, 0) node (g2_1){};
    \path (g2) +(0, 3/2 * \size) node (g2_2){};
    \path (g2) +(\size, 5/2 * \size) node (g2_3){};
    \path (g2) +(2 * \size, 3/2 * \size) node (g2_4){};
    \path (g2) +(2 * \size, 0) node (g2_5){};
    \path (g2) +(1/2 * \size, 1/2 * \size) node (g2_6){};
    \path (g2) +(1/2 * \size, \size) node (g2_7){};
    \path (g2) +(\size, 3/2 * \size) node (g2_8){};
    \path (g2) +(3/2 * \size, \size) node (g2_9){};
    \path (g2) +(3/2 * \size, 1/2 * \size) node (g2_10){};

    \foreach \Point in {(g2_1), (g2_2), (g2_3), (g2_4), (g2_5)}{
        \node at \Point{\textcolor{blue}{\textbullet}};
    }
    \foreach \Point in {(g2_6), (g2_7), (g2_8), (g2_9), (g2_10)}{
        \node at \Point{\textcolor{red}{\textbullet}};
    }
    \draw [blue]
            (g2_1) -- (g2_2)
            (g2_1) -- (g2_5)
            (g2_2) -- (g2_3)
            (g2_3) -- (g2_4)
            (g2_4) -- (g2_5);
    \draw [red]
            (g2_6) -- (g2_8)
            (g2_6) -- (g2_9)
            (g2_7) -- (g2_9)
            (g2_7) -- (g2_10)
            (g2_8) -- (g2_10);
    \draw [blue]
            (g2_1) -- (g2_6)
            (g2_2) -- (g2_7)
            (g2_3) -- (g2_8)
            (g2_4) -- (g2_9)
            (g2_5) -- (g2_10);
    \path (g2) ++(\size,-\size / 2) node[draw=none,fill=none] { {\large (b)}};
}

\end{tikzpicture}
\end{center}